\documentclass[12pt]{amsart}
\newcounter{defcounter}
\usepackage{amssymb,amsmath,amscd,graphicx, enumitem,
latexsym,amsthm}
\usepackage{amssymb,latexsym,eufrak,amsmath,amscd,graphicx}
  \usepackage[all]{xy}
  \usepackage{pdfsync}
\theoremstyle{plain}
\newtheorem{theorem}{Theorem}
\newtheorem*{reciprocity}{Hermite Reciprocity}
\newtheorem{proposition}[theorem]{Proposition}
\newtheorem{corollary}[theorem]{Corollary}

\newtheorem{proposition.definition}[theorem]{Proposition/Definition}

 \newtheorem{folkconjecture}[theorem]{Folk-Conjecture}

\newtheorem{conjecture}[theorem]{Conjecture}

\theoremstyle{definition}
\newtheorem{definition}[theorem]{Definition}

\newtheorem{remark}[theorem]{Remark}
\newtheorem{example}[theorem]{Example}

\newcommand{\lra}{\longrightarrow}

\newcommand{\noi}{\noindent}
\newcommand{\PP}{\mathbf{P}}

\newcommand{\CC}{\mathbf{C}}

\newcommand{\OO}{\mathcal{O}}

\newcommand{\FF}{\mathcal{F}}

\newcommand{\rk} {\text{rank }}

\newcommand{\GGG}[2]{\Gamma \big( {#1},{#2}\big)}

\newcommand{\pr}{\prime}

\newcommand{\Hilb}{\textnormal{Hilb}}

\newcommand{\Cliff}{\textnormal{Cliff}}
\newcommand{\Moduli}{\mathfrak{M}}
\newcommand{\Tan}{\textnormal{Tan}}
\newcommand{\SL}{\textnormal{SL}}

\newcommand{\AAA}{\textbf{A}}

\numberwithin{theorem}{section}
\numberwithin{equation}{section}

\begin{document}

\title[Tangent Developable Surfaces and Defining Equations]
{Tangent Developable Surfaces and the Equations Defining Algebraic Curves}
 \author{Lawrence Ein}
 \address{Department of Mathematics, University of Illinois at Chicago, 851 South Morgan St., Chicago, IL  60607}
  \email{{\tt ein@uic.edu}}
  \thanks{Research of the first author partially supported by NSF grant DMS-1801870.}
 
 \author{Robert Lazarsfeld}
  \address{Department of Mathematics, Stony Brook University, Stony Brook, New York 11794}
 \email{{\tt robert.lazarsfeld@stonybrook.edu}}
 \thanks{Research  of the second author partially supported by NSF grant DMS-1739285.}

\maketitle

\begin{abstract} This  is  an introduction, aimed at a general mathematical audience,  to  recent work of  Aprodu, Farkas, Papadima, Raicu and Weyman. These authors established  a long-standing folk conjecture concerning the equations defining the tangent developable surface of the rational normal curve. This in turn led to a new proof of a fundamental theorem of Voisin concerning the syzygies of generic canonical curves. The present note, which is the   write-up of a talk given by the second author at the Current Events seminar at the 2019 JMM, surveys this circle of ideas. \end{abstract}
    
 \section*{Introduction}

Let $X$ be a smooth complex projective algebraic curve -- or equivalently a compact Riemann surface -- of genus $g \ge 2$, and denote by
\[   H^{1,0}(X) \ = \ \GGG{X}{\Omega_X^1} \]
the $\CC$-vector space of holomorphic $1$-forms on $X$. Recalling that $\dim H^{1,0}(X) = g$, choose a basis
\[   \omega_1 \ldots , \omega_{g} \, \in \, H^{1,0}(X). \]
It is classical that the $\omega_i$ do not simultaneously vanish at any point $x\in X$, so one can define a holomorphic map
\[   \phi_X : X \lra \PP^{g-1} \ \ ,  \ \ x \mapsto [\, \omega_1(x) , \ldots , \omega_{g}(x)\,  ] \]
from $X$ to a projective space of dimension $g -1$, called the \textit{canonical mapping} of $X$.  With one well-understood class of exceptions, $\phi_X$ is an embedding, realizing $X$ as an algebraic curve 
\[  X \ \subseteq \ \PP^{g-1} \]
 of degree $2g-2$. Any compact Riemann surface admits many  projective embeddings, but the realization  just constructed has the signal advantage of being  canonically defined up to a linear change of coordinates on $\PP^{g-1}$. Therefore  the extrinsic projective geometry of a canonically embedded curve must reflect its intrinsic geometry, and working this principle out is an important theme in the theory of algebraic curves. 
 
 Given any projective variety, one can consider the degrees of its defining equations. An important   theorem of Petri from 1922 states that with a slightly wider range of exceptions, a canonical curve $X \subseteq \PP^{g-1}$ is cut out by quadrics, i.e. polynomials of degree two. Classically that seemed to be the end of the story, but in the early 1980s  Mark Green realized that Petri's result should be the first case   of a much more general statement  involving higher syzygies. In other words, one should consider not only the defining equations themselves, but  the relations among them,  the  relations among the relations, and so on. The resulting conjecture has attracted a huge amount of attention over the past thirty-five years.
 
As of this writing, Green's conjecture remains open. However Voisin made a major breakthrough in 2002 by proving that it holds for \textit{general} curves, where one rules out for instance all the sorts of exceptional cases alluded to above. Her proof introduced a number of very interesting new ideas, but at the end of the day it relied on some difficult and lengthy cohomological calculations. Prior to Voisin's work, O'Grady and  Buchweitz--Schreyer had observed that one might be able to attack the syzygies of generic canonical curves by studying a very concrete and classical object, namely the developable surface of tangent lines to a rational normal curve. A substantial body of experimental evidence supported this proposal, but in spite of considerable effort  nobody was able to push through the required computations. In a recent preprint \cite{AFPRW}, however,  Aprodu, Farkas, Papadima, Raicu and Weyman have succeeded in doing so. Their work is the subject of the present report.  

This note is organized as follows. Section 1 is devoted to the geometry of canonical curves, the basic ideas around syzygies, and the statement of Green's conjecture. The case of general  curves, and its relation to the tangent surface of rational normal curves occupies \S 2. Finally, in \S3 we explain the main  ideas underlying the work of AFPRW from an algebro-geometric perspective.   

We work throughout over the complex numbers. In particular, we completely ignore contributions of \cite{AFPRW} to understanding what parts of Green's conjecture work in positive characteristics.

We thank the authors of \cite{AFPRW} for sharing an early draft of their paper. We have profited from correspondence and conversations  with David Eisenbud, Gabi Farkas, Claudiu Raicu, Frank Schreyer and Claire Voisin.

\section{Canonical curves, syzygies and Green's conjecture}

\subsection*{Canonical curves and Petri's theorem} Denote by $X$ a smooth complex projective curve of genus $g\ge 2$, and as in the Introduction consider the canonical mapping
\[   \phi_X  \, : \, X \lra \PP^{g-1} \ \ , \ \ x \mapsto [ \, \omega_1(x), \ldots, \omega_g(x) ], \]
arising from a basis $\omega_1, \ldots , \omega_g \in H^{1,0}(X)$  of holomorphic  $1$-forms on $X$. By construction the inverse image of a hyperplane cuts out the zero-locus of a such a $1$-form, and therefore consists of $2g-2$ points (counting multiplicities).  It is instructive to consider concretely the first few cases.

\begin{example} \textbf{(Genus $2$ and hyperelliptic curves)}. Suppose $g(X) = 2$. Then the canonical mapping is a degree two branched covering
\[ \phi_X \, : \, X \lra \PP^1. \]
In general, a curve of genus $g$ admitting a degree two covering $X \lra \PP^1$ is called \textit{hyperelliptic}. Thus every curve of genus $2$ is hyperelliptic, but when $g \ge 3$ these are in many respects the ``most special" curves of genus $g$. The canonical mapping of a hyperelliptic curve factors as the composition
\[   X \, \lra \, \PP^1 \, \subseteq  \, \PP^{g-1} \]
of the hyperelliptic involution with an embedding of $\PP^1$ into $\PP^{g-1}$. It follows from the Riemann--Roch theorem that these are the only curves for which $\phi_X$ is not an embedding:

\begin{center}\parbox{5in} {\textsc{Fact.} If $X$ is non-hyperelliptic, then the canonical mapping 
\[ \hskip -2in\phi_X : X \ \subseteq \ \PP^{g-1} \] is an embedding.}
\end{center}
\end{example}

\begin{example}\textbf{(Genus 3 and 4).} Assume that $X$ is not hyperelliptic. 
When $g(X) = 3$, the canonical mapping realizes $X$ as a smooth curve of degree $4$ in $\PP^2$, and any such curve is canonically embedded. When $g = 4$, $\phi_X$ defines an embedding
\[   X \, \subseteq \, \PP^3  \]
in which $X$ is the complete intersection of a surface of degree $2$ and degree $3$.
\end{example}

\begin{example} \textbf{(Genus 5).} This is the  first case where one sees the interesting behavior of quadrics through a canonical curve. Consider a canonically embedded non-hyperelliptic curve of genus $5$
\[  X \ \subseteq \ \PP^4 \ \ ,  \ \deg X = 8. \]
One can show that there is a three-dimensional vector space of quadrics through $X$, say with basis $Q_1, Q_2, Q_3$. There are now two possibilities:
\begin{itemize}
\item[(a).] $X$ is \textit{trigonal}, i.e. there exists a degree three branched covering
\[ \pi : X \lra \PP^1. \]
In this case each of the fibres of $\pi$ spans a line in $\PP^{4}$, and hence any quadric containing $X$ must also contain each of these lines. They sweep out a ruled surface $S \subseteq \PP^4$ containing $X$, and three quadrics through $X$ meet precisely along $S$: \[   Q_1 \, \cap \, Q_2 \, \cap \, Q_3 \ = \ S.\]
The canonical curve  $X$ is cut out in $S$ by some cubic forms. 
\vskip 10pt
\item[(b).] $X$ is not trigonal, i.e. cannot be expressed as a $3$-sheeted branched covering of $\PP^1$. Then $X$ is the complete intersection of the the three quadrics containing it:
\[  X \ = \ Q_1 \, \cap \, Q_2 \, \cap \, Q_3. \]
This is the general case.
\end{itemize}
\end{example}

\begin{example} \textbf{(Genus 6).} Consider finally a non-hyperelliptic canonical curve $X \subseteq \PP^5$ of genus $6$. Now the polynomials of degree two vanishing on $X$ form a vector space of dimension $6$, and there are three cases:
\begin{itemize}
\item[(a).] If $X$ is {trigonal}, then as above the quadrics through $X$ intersect along  the ruled surface $S$ swept out by the trigonal divisors.
\vskip 10pt
\item[(b).] Suppose that $X$ is a smooth curve of degree $5$ in $\PP^2$. In this case the canonical image of $X$ lies on the \textit{Veronese surface} 
\[   X \ \subseteq \ V \ \subseteq \ \PP^5, \]
a surface of degree $4$ abstractly isomorphic to $\PP^2$, and $V$ is the intersection of the quadrics through $X$ in canonical space.
\vskip 10pt
\item[(c).] The general situation is that $X$ is neither trigonal nor a plane quintic, and then $X \subseteq \PP^5$ is cut out by the quadrics passing through it. Note however that $X$ is not the complete intersection of these quadrics, since they span a vector space of dimension  strictly greater than the codimension of $X$.
\end{itemize}
\end{example}

We conclude this subsection by stating Petri's theorem. Consider a non-hyperelliptic canonical curve $X \subseteq \PP^{g-1}$. Let $S = \CC[ Z_0, \ldots, Z_{g-1} ]$ be the homogeneous coordinate ring of  canonical space $\PP^{g-1}$, and denote by  
\[  I_X \ \subseteq \  S\]
the homogeneous ideal of all forms vanishing on $X$. We ask when $I_X$ is generated by forms of degree $2$: this is the strongest sense in which $X$ might be cut out by quadrics.
\begin{theorem} [Petri]  
The homogeneous ideal $I_X$ fails to be generated by quadrics if and only if $X$ is either trigonal or a smooth plane quintic. 
\end{theorem}
\noi Note that the Petri-exceptional curves fall into two classes: there is one family (trigonal curves) that appears in all genera, and in addition one ``sporadic" case.

In retrospect, Petri's statement suggests some natural questions. For example, how does one detect algebraically curves $X$ that can be expressed as  a degree $4$ branched covering $X \lra \PP^1$, or that arise as smooth plane sextics? Or again, what happens in the generic case, when $X$ does not admit any unusually low degree mappings to projective space? Green's beautiful insight is that one should consider for this not just the generators of $I_X$ but also its higher syzygies.

\subsection*{Syzygies} The idea to study  the relations -- or syzygies -- among the generators of an ideal goes back to  Hilbert.  Making this  precise  inevitably involves a certain amount of notation, so  perhaps it's best to start concretely with the simplest example.\footnote{We refer the reader to \cite{EisenbudBook} for an systematic introduction to the theory from an algebraic perspective, and \cite{Aprodu.Nagel} for a more geometric viewpoint.}

The rational normal curve $C \subseteq \PP^3$ of degree $3$ is the image of the embedding
\[  \nu : \PP^1 \lra \PP^3 \ \ , \ \ [ \, u , v \,] \mapsto [ \, u^3, u^2v, uv^2 , v^3 ].\]
Writing $[\, Z_0, \ldots, Z_3 \, ]$ for homogeneous coordinates on $\PP^3$, it is  a pleasant exercise to show that $C$ can be described as the locus where a catalecticant matrix drops rank:
\[
C \ = \  \left \{ \ \rk \! \! \! \left [ \begin{matrix}   Z_0 &Z_1 &Z_2 \\ Z_1 & Z_2  &Z_3  \end{matrix} \right ]\le 1 \ \right \}.
\]
Therefore $C$ lies on the three quadrics
\[   Q_{02} \, = \, Z_0Z_2- Z_1^2 \ , \  Q_{03} \, = \, Z_0Z_3 - Z_1Z_2 \ , \ Q_{13} \, = \, Z_1Z_3 - Z_2^2 \]
given by the $2 \times 2$ minors of this matrix, and in fact these generate the homogeneous ideal $I_C $ of $C$. While the $Q_{ij}$ are linearly independent over $\CC$, they satisfy two relations with polynomial coefficients, namely\[  \begin{aligned}   Z_0 \cdot Q_{13} - Z_1 \cdot Q_{03}  + Z_2 \cdot Q_{02} \, &= \, 0 \\ Z_1 \cdot Q_{13} - Z_2 \cdot Q_{03}  + Z_3 \cdot Q_{02} \, &= \, 0 . 
\end{aligned}\tag{*}
\]
One can derive these by repeating a row of the matrix defining $C$ and expanding the resulting determinant along the duplicate row. Moreover it turns out that any relation among the $Q_{ij}$ is a consequence of these. 

We  recast this discussion somewhat more formally. Write $S = \CC[Z_0, \ldots, Z_3]$ for the homogeneous coordinate ring of $\PP^3$. The three quadric generators of $I_C$ determine a surjective map
\[   
 S(-2)^{\oplus 3}  \lra I_C,
\]
where $S(-2)$ denotes a copy of $S$ re-graded so that multiplication by the $Q_{ij}$ is degree preserving. The relations in (*) come from choosing generators for the kernel of this map. So the upshot of the previous paragraph is that one has an exact sequence
\[
\xymatrix{ 0 \ar[r] &S(-3)^{\oplus 2} \ar[rr] ^{\left(\begin{smallmatrix}{ \ }Z_0 & {\ }Z_1\\-Z_1 &{  }	 -Z_2\\  {\ }Z_2 & {\ }Z_3\end{smallmatrix}\right)}& & S(-2)^{\oplus 3} \ar[rrr] ^{\ \ \left(\begin{smallmatrix}Q_{13} \ & Q_{03} \ & Q_{02} \end{smallmatrix}\right)}& & & I_C \ar[r] &0 
}
\]
of $S$-modules. This is the \textit{minimal graded free resolution} of $I_C$. 

The general situation is similar. Sticking for simplicity to the one-dimensional case, consider a non-degenerate curve 
\[  C \ \subseteq  \ \PP^r\] i.e.\! one not lying on any hyperplanes. We suppose in addition that $C$ is projectively normal, a  technical condition that holds for any embedding of sufficiently large degree (and for non-hyperelliptic canonical curves thanks to a theorem of Noether.) Put \[ S \ = \ \CC[Z_0, \ldots, Z_r],\]and denote by $I_C \subseteq S$ the homogeneous ideal of $C$. Then $I_C$ has a minimal resolution $E_{\bullet}$ of length $r-1$:
\begin{equation}\label{Resolution}
 \xymatrix{
0 \ar[r] & E_{r-1} \ar[r]& \ldots \ar[r]  & E_2  \ar[r]  & \ar[r] E_1 \ar[r]  & I_C  \ar[r]&0 ,}
 \end{equation}
where
$  E_i = \oplus \, S(-a_{i,j})$.  It is elementary that $a_{i,j} \ge i+1$ for every $j$. 

Green realized that the way to generalize classical statements about quadratic generation of $I_C$ is to ask when the first $p$ steps of this resolution are as simple as possible.
\begin{definition}
One says that $C$ satisfies Property ($N_p$) if 
\[   E_i \ = \ \oplus \,  S(-i-1) \ \ \text{ for every } \ 1 \le i \le p. \]
\end{definition}
\noi Thus ($N_1$) holds if and only if $I_C$ is generated in degree $2$.  The first non-classical condition is ($N_2$), which asks in addition that if one chooses quadratic generators $Q_\alpha \in I_X$, then the module of syzygies among the $Q_\alpha$ should be spanned by relations of the form
\[  \sum \, L_\alpha \cdot Q_\alpha \ = \ 0,  \tag{*} \]
where the $L_\alpha$ are \textit{linear} polynomials. Condition ($N_3$) would ask that the syzygies among the coefficient vectors describing the relations (*)  are themselves generated by polynomials of degree one. 

\begin{example}  The twisted cubic $C \subseteq \PP^3$ discussed above satisfies ($N_2$). On the other hand, an elliptic curve $E \subseteq \PP^3$ of degree $4$  is the complete intersection of two quadrics, whose homogeneous ideal has a Koszul resolution:
\[  0 \lra S(-4) \lra S(-2)^{\oplus 2} \lra I_E \lra 0. \]
Thus $E$ satisfies ($N_1$) but not ($N_2$). 
\end{example}

Return now to a non-hyperelliptic canonical curve $X \subseteq \PP^{g-1}$ of genus $g$.  Petri's theorem states that $X$ satisfies ($N_1$) unless it is trigonal or a smooth plane quintic. Green's conjecture vastly extends this by predicting when $X$ satisfies condition ($N_p$). 

\subsection*{Green's conjecture} In order to state Green's conjecture, it remains  to understand the  pattern behind the exceptional cases in Petri's theorem.

 Let $X$ be a curve of genus $g \ge 2$, and suppose given a non-constant holomorphic mapping
\[   \varphi: X \lra \PP^r. \]
We assume that $X$ does not map into any hyperplanes, in which case we write $r(\varphi) = r$: this is often called the \textit{dimension} or \textit{rank} of $\varphi$. If $\varphi$ has degree $d$ in the sense that a general hyperplane pulls back to $d$ points on $X$, we set $d(\varphi) = d$. The \textit{Clifford index}  of $\varphi$ is then defined to be
\[  \Cliff(\varphi) \ = \ d(\varphi) - 2 \cdot r(\varphi). \]
A classical theorem of Clifford states that if $d(\varphi) \le g-1$, then 
\[  \Cliff(\varphi) \, \ge \, 0, \]
and equality holds if and only if $X$ is hyperelliptic and $\varphi : X \lra \PP^1$ is the hyperelliptic involution (or a mapping derived from it by a Veronese-type construction). 

We now attach an invariant to $X$ by considering the minimum of the Clifford indices of all ``interesting" mappings:
\begin{definition} The Clifford index of $X$ is
\[ \Cliff(X) \ = \ \min \big \{ \Cliff(\varphi) \mid \, d(\varphi) \le g-1 \, \big \}.
\]
\end{definition}
\noi One has
\[  0 \ \le \ \Cliff(X) \ \le \  \left [ \frac{g-1}{2} \right ], \]
 for every $X$, the first inequality coming from Clifford's theorem, and the second (as we explain  in the next section) from Brill-Noether theory. 
 Moreover
\[   \Cliff(X) \, = \, 0 \ \ \Longleftrightarrow \ \ X \text{ is hyperelliptic}, \]
and similarly one can show that when $X$ is non-hyperelliptic,
\[  
\Cliff(X) \, = \, 1 \ \ \Longleftrightarrow \ \ X \text{ is either trigonal or a smooth plane quintic}. \]

It is now clear what to expect for higher syzygies:
\begin{conjecture} [Green, \cite{Green}] \label{Green}
Let $X \subseteq \PP^{g-1}$ be a non-hyperelliptic canonical curve. Then the Clifford index of $X$ is equal to the least integer $p$ for which Property $(N_p)$ fails for $X$.
\end{conjecture}
\noi  The case $p = 1$ is exactly Petri's theorem, and the first non-classical case $p= 2$ was established by Schreyer \cite{Schreyer} and Voisin \cite{Voisin0}. There is a symmetry among the syzygies of canonical curves, and knowing the smallest value of $p$ for which ($N_p$) fails turns out to determine the grading of the whole resolution of $I_X$.

One implication in Green's statement is elementary: it was established in the appendix to \cite{Green} that if $\Cliff(X) = p$, then ($N_p$) fails for $X$. 
What remains mysterious as of this writing is how to show conversely that unexpected syzygies actually have a geometric origin. 

\begin{remark} [The gonality conjecture] Inspired by Conjecture \ref{Green}, Green and the second author proposed in \cite{GreenLaz} that one should be able to read off the gonality of a curve from the resolution of any one line bundle of sufficiently large degree. This was originally envisioned as a warm-up problem for canonical curves, but the present authors observed a few years ago that in fact a small variant of the Hilbert-schematic ideas introduced by Voisin in \cite{Voisin1} leads to a very quick proof of the conjecture \cite{Gonality}. See \cite{SyzSurvey} for a survey of recent work about syzygies of varieties embedded by very positive linear series.
\end{remark}

\section{General  curves of large genus and the tangent developable to rational normal curves}

\subsection*{General curves}  The most important instance of Green's conjecture -- which is the actual subject of the present report -- is the case of ``general" curves. We start by explaining a little more precisely what one means by this.

In the 1960s, Mumford and others constructed an algebraic variety $\Moduli_g$ whose points parameterize in a natural way isomorphism classes of smooth projective curves of genus $g \ge 2$. This is the \textit{moduli space} of curves of genus $g$. One has \[ \dim \, \Moduli_g \ = \ 3g - 3, \] formalizing a computation going back to Riemann that compact Riemann surfaces of genus $g\ge 2$ depend on $3g - 3$ parameters. Special classes of curves correspond to   (locally closed) proper subvarieties of $\Moduli_g$: for example, hyperelliptic curves are parameterized by a subvariety $\mathfrak{H}_g \subseteq \Moduli_g$ of dimension $2g -1$, showing again that hyperelliptic curves are special when $g \ge 3$.  
One says that a statement holds for a \textit{general  curve} of genus $g$ if it holds for all curves whose moduli points lie outside a finite union of proper subvarieties of $\Moduli_g$.

The question of what mappings $\varphi : X \lra \PP^r$ exist for a general curve $X$ was studied classically, and the theory was put on a firm modern footing in the 1970s by Kempf, Kleiman-Laksov, Griffiths-Harris and Gieseker, among others: see \cite{ACGH}. For our purposes, the basic fact is the following:
\begin{theorem} [Weak form of Brill--Noether theorem] \label{BN.Thm} Let $X$ be a general curve of genus $g \ge 2$. Then there exists a map
$ \varphi : X \lra \PP^r $
of degree $d$ and dimension $r$ if and only if
\[   g \ \ge \ (r+1)(g - d + r). \]
In particular,
$\Cliff(X) = \left[ \frac{g-1}{2} \right]. $
\end{theorem}

Green's conjecture then predicts the shape of the minimal resolution of the ideal of a general canonical curve of genus $g$. This is the stunning result established by Voisin \cite{Voisin1}, \cite{Voisin2}.
\begin{theorem} [Voisin's Theorem] \label{Voisin.Thm} Put $c = [ \frac{g-3}{2}].$ Then a general canonical curve $X \subseteq \PP^{g-1}$ of genus $g$ satisfies Property $(N_c)$.
\end{theorem}
\noi  The symmetry mentioned following the statement of Green's conjecture imposes limits on how far ($N_p$) could be satisfied, and one  can view Voisin's theorem as asserting that the syzygies of a general canonical curve are ``as linear as possible" given this constraint. 

General principles imply that the set of curves for which the conclusions of Theorems \ref{BN.Thm} or \ref{Voisin.Thm} hold are parameterized by Zariski-open subsets of $\Moduli_g$. So to prove the results  it would suffice to exhibit one curve of each genus $g$ for which the assertions are satisfied. However it has long been understood that this is not a practical approach.
Instead, two different strategies have emerged for establishing statements concerning general canonical curves. 

The first is to consider singular rational curves. For example, a rational curve $\Gamma \subseteq \PP^{g-1}$ of degree $2g - 2$ with $g$ nodes can be realized as a limit of canonical curves. The first proof  of Theorem \ref{BN.Thm}, by Griffiths and Harris,  went by establishing that $\Gamma$ satisfies an appropriate analogue of the statement, and then deducing that \ref{BN.Thm} must hold for a general smooth curve of genus $g$. A difficulty here is that the nodes themselves have to be in general position, requiring a further degeneration. Eisenbud and Harris subsequently found that it is much  better to work with cuspidal curves: we will return to this shortly. More recently, tropical methods have entered the picture to give new proofs of Theorem \ref{BN.Thm}.

A different approach, initiated  in \cite{BNP}, involves  K3 surfaces. These are surfaces \[ S =S_{2g-2} 
\ \subseteq \ \PP^g\] of degree $2g-2$ whose hyperplane sections are canonical curves. It turns out to be quite quick to show that  these curves are Brill-Noether general provided that $S$ itself is generic. While it is not easy to exhibit explicitly a suitable K3, it is known by Hodge theory that they exist in all genera. This re-establishes the existence of curves that behave generically from the perspective of Brill-Noether theory.  

This was the starting point of Voisin's proof of Theorem \ref{Voisin.Thm}. Under favorable circumstances the resolution of a surface restricts to that of its hyperplane section, so it suffices to show that a general K3 surface of genus $g$ satisfies the conclusion of \ref{Voisin.Thm}. However so far this doesn't really simplify the picture. Voisin's remarkable new idea was to pass to a larger space, namely the Hilbert scheme 
\[ S^{[c+1]} \ = \ \Hilb^{c+1}(S) \]  parameterizing finite subschemes of length $(c+1)$ on $S$. Voisin showed that the syzygies of $S$ are encoded in a quite simple-looking geometric statement on $S^{[c+1]}$. The required computations turn out to be  rather involved, but 
in a real tour de force Voisin succeeded in pushing them through. Interestingly, it later emerged that her computations could be used to establish many other cases of Green's conjecture, eg that it holds for a general curve of each gonality, or for every curve appearing on a K3 surface. See \cite{Aprodu.Nagel} for some examples and references.

At about the same time that Green formulated his conjecture in the early 1980s, Eisenbud and Harris \cite{EH} realized that many of the difficulties involved in degenerating to nodal rational curves  disappeared if one worked instead with rational curves 
 with $g$ cusps.  The advantage of these curves is that they behave Brill-Noether generally without any conditions on the location of the singular points. This raised the possibility that one might use $g$-cuspidal curves to study syzygies of general canonical curves. It was at this point that Kieran O'Grady, and independently Buchweitz and Schreyer, remarked  that it should suffice to understand the syzygies of a very classical object, namely the tangent surface to a rational normal curve.

\subsection*{The tangent developable of a rational normal curve} 
\begin{figure}
\includegraphics[width=.8\hsize]{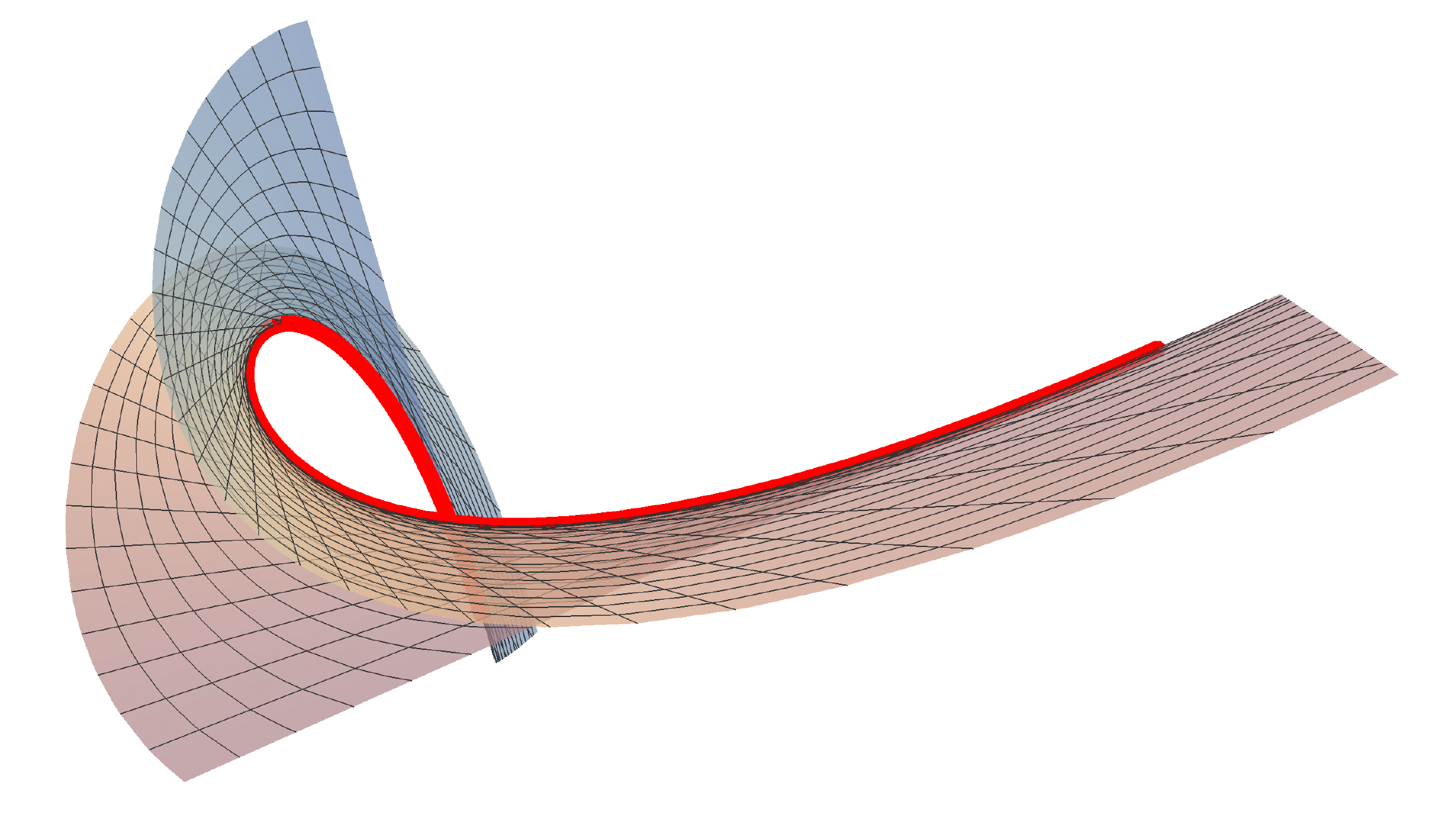}
\caption{The tangent developable surface to the twisted cubic in $\PP^3$}
\label{Fig1}
\end{figure}

Let $C \subseteq \PP^g$ be a rational normal curve of genus $g$. By definition this is the image of the embedding  
\[  \PP^1 \hookrightarrow \PP^g \ \ \text{ given by } \ \  [s,t] \, \mapsto \, [\, s^g\, , \, s^{g-1}t\, , \, \ldots \, , \, st^{g-1} \, , \, t^g \, ] . \]
One can associate to $C$ (as to any smooth curve) its \textit{tangent surface}
\[   T \ = \ \textnormal{Tan}(C) \ \subseteq \ \PP^g, \]
defined to be the union of all the embedded projective tangent lines to $C$. In the case at hand, one can describe $T$ very concretely. Specifically it is the image of the map 
\begin{equation} \label{Normalization}
 \nu \, :  \,  \PP^1 \times \PP^1 \lra \PP^g 
 \end{equation} given matricially by 
\begin{gather} \label{Tang.Develop.Formula}
\nu  \big( [s,t] \times [u , v] \big) \ = \   \left[\, u \ \ v \,  \right] \cdot \textnormal{Jac}(\mu), \end{gather}
where $\textnormal{Jac}(\mu)$ is the $2 \times (g+1)$ matrix of partials of $\mu = [ \, s^g, s^{g-1}t, \ldots ,st^{g-1},  t^g]$. In other words, 
\begin{equation*}  \nu  \big( [s,t] \times [u,v] \big) \ = \ \Big[ \, g\cdot 
s^{g-1}u \  , \ (g-1)\cdot s^{g-2}tu + s^{g-1}v \ , \ \ldots \ , \ g\cdot st^{g-1}v \, \Big] .
\end{equation*}
Note that $\nu$ is one-to-one, and maps the diagonal $\Delta \subseteq \PP^1 \times \PP^1$ isomorphically to $C$. However $\nu$ is not an embedding: it ramifies along the diagonal,  and $T$ has cuspidal singularities along $C$. Figure \ref{Fig1} shows a drawing of $T$ in the case $g = 3$.  
The tangent surface $T$ is a  complex-geometric analogue of one of the classes of developable surfaces studied in differential geometry. A pleasant computation shows that  $\deg(T) = 2g -2$.\footnote{Either observe that $\nu$ is given by an (incomplete)  linear series of type $(g-1,1)$ on $\PP^1 \times \PP^1$, or use Riemann-Hurwitz for a degree $g$ mapping $\PP^1 \lra \PP^1$ to see that there are $2g-2$ tangent lines to $C$ meeting a general linear space $\Lambda \subseteq \PP^g$ of codimension $2$.} 

The upshot of this discussion is that the hyperplane sections of $T$ are rational curves $\Gamma\subseteq \PP^{g-1}$ of degree $2g-2$ with $g$ cusps -- in other words, the degenerations of canonical curves with which one hopes to be able to prove the generic case of Green's conjecture. This led to the
\begin{folkconjecture} \label{Folk.Conj}
The tangent developable surface 
\[ T \  = \ \textnormal{Tan}(C) \ \subseteq \ \PP^g \] satisfies Property $(N_p)$ for 
$p= \left[ \frac{g-3}{2} \right ]$.
\end{folkconjecture}
\noi 
With a small argument showing that $T$ indeed has the same syzygies as its hyperplane sections, it has been well understood since the mid 1980s that this would imply the result (Theorem \ref{Voisin.Thm}) that Voisin later proved by completely different methods. 

The important thing to observe about \ref{Folk.Conj} is that it is   a completely concrete statement. Via the parameterization \eqref{Tang.Develop.Formula}, the conjecture was  quickly verified   for a large range of genera  using early versions of the computer algebra system Macaulay.  That such an utterly down-to-earth assertion could resist proof for thirty-five years has been something of a scandal. Happily, the work of  Aprodu, Farkas, Papadima, Raicu and Weyman has remedied this situation.

\section{Sketch of the proof of Conjecture \ref{Folk.Conj} }

In this section, we outline the main ideas of the work of AFPRW proving Folk Conjecture \ref{Folk.Conj}. 

The actual write-up in \cite{AFPRW} is a bit long and complicated, in part because the authors work to extend their results as far as possible to positive characteristics, and in part because they are fastidious in checking that the maps that come up are the expected ones. Here we focus  on the essential geometric  ideas that seem to underlie their computations.

\subsection*{Computing the syzygies of $T$}

The first step in the argument of \cite{AFPRW} is to understand the tangent developable $T = \Tan(C)$ and its syzygies in terms of more familiar and computable objects. This culminates in Theorem \ref{Syzygy.Tang.Summary} below, which describes the relevant syzygies linear algebraically. Some of the computations of AFPRW apparently elaborate on earlier (unpublished) work of Weyman, as outlined in Eisenbud's notes \cite{Eisenbud}.

A basic principle guiding algebraic geometry holds that spaces are determined by the polynomial functions on them, so we will need to understand those on $T$.  It is in turn natural to expect that functions on the tangent developable should be described using the mapping  $\nu  :  \PP^1 \times \PP^1 \lra T$ from \eqref{Normalization}, which realizes $T$ as the homeomorphic image of $\PP^1 \times \PP^1$ cusped along the diagonal. In order to get a sense of how this should go, let us start with a one-dimensional toy example. 

\begin{figure}
\includegraphics[scale = .45]{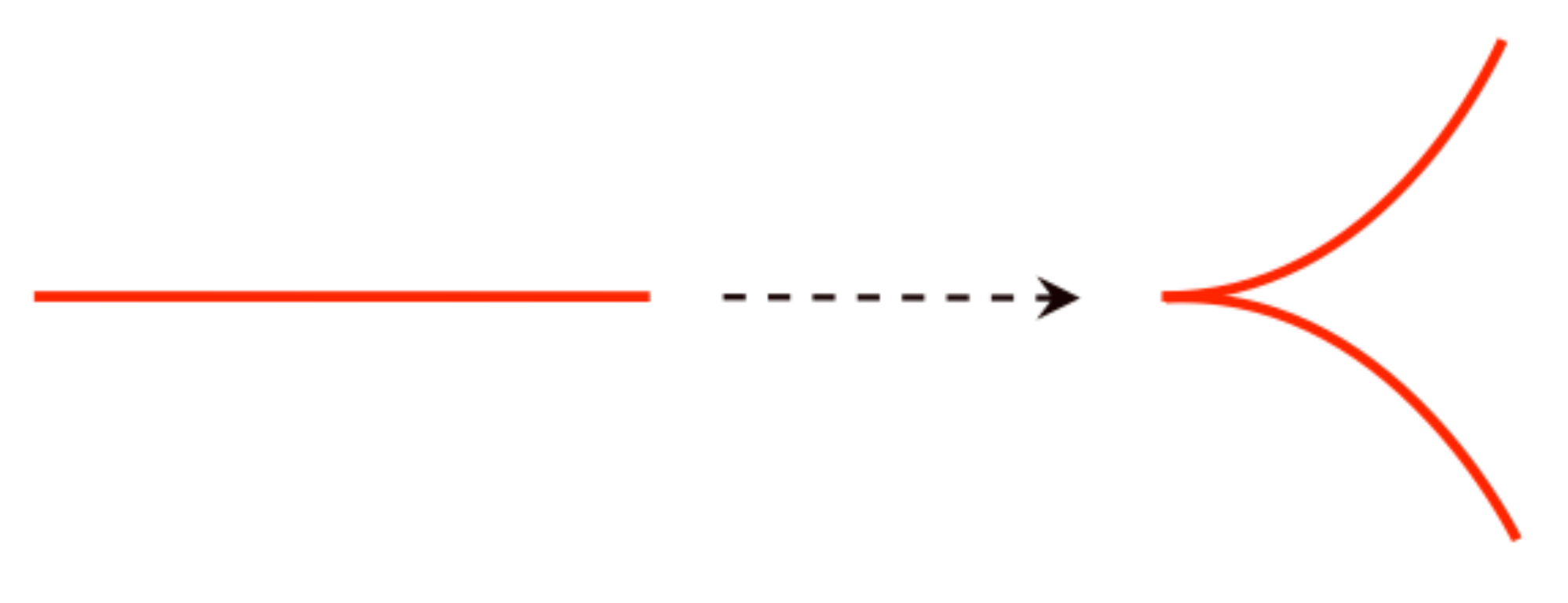}
\caption{The mapping $\nu_0(t) =  ( t^2 \, , \, t^3 )$}
\label{Fig2}
\end{figure}

Consider then the mapping
\[  \nu_0 : \AAA^1 = \CC  \lra \AAA^2 = \CC^2\ \ , \ \ \nu_0(t) = \big( t^2 \, , \, t^3 \big). \]
This maps $\AAA^1$ homeomorphically onto the cuspidal curve 
\[   T_0 \, = \,   \{ y^2 \, = \, x^3   \} \ \subseteq \ \AAA^2 \]
in the plane, and the polynomial functions on $T_0$ are realized as the subring
\[  \CC[T_0] = \CC[ t^2, t^3 ] \ \subseteq \ \CC[t] = \CC[\AAA^1]\]of the regular functions on the affine line. (See Figure \ref{Fig2}.) 

The point to note is that we can describe $\CC[T_0]$  intrinsically, without using the map $\nu_0$. Specifically, there is a $\CC$-linear derivation
\[   \delta_0 \, : \, \CC[t] \lra \CC \ \ , \ \ \delta_0(f) = f^\pr (0), \]
and $\CC[T_0] = \ker(\delta_0)$. Moreover while $\delta_0$ is not $\CC[t]$-linear, it is linear over $\CC[t^2, t^3]$, giving a short exact sequence 
$0 \lra \CC[T_0] \lra \CC[\AAA^1] \lra \CC \lra 0$ of $\CC[T_0]$-modules.

This model  generalizes. Writing $\OO_X$ to denote the (sheaf of locally) polynomial  functions on a variety $X$, one has
\begin{proposition} \label{OT.as.Kernel}
There is a $\CC$-linear derivation
\[ \delta : \OO_{\PP^1 \times \PP^1} \lra \Omega^1_\Delta \ \ , \ \  \delta(f) = df | \Delta
\]
with $\ker \delta = \OO_T$. Moreover, this gives rise to a short exact sequence
\begin{equation} \label{Basic.Exact.Seq}
 0 \lra \OO_T \lra \OO_{\PP^1 \times \PP^1} \overset{\delta} \lra \Omega^1_C \lra 0 \end{equation}
$($of sheaves$)$ on $\PP^g$.\footnote{Strictly speaking, the middle term in (*) is  the direct image  $\nu_* \big(\OO_{\PP^1 \times \PP^1}\big)$, but we wish to minimize sheaf-theoretic notation.}
\end{proposition}
\noi It is easy to describe the syzygies of $\OO_{\PP^1 \times \PP^1}$ and $\Omega^1_C$, and then the plan is  to use \eqref{Basic.Exact.Seq} to analyze the syzygies of $T$. 

At this point we require some additional syzygetic notation. As above denote by $S = \CC[Z_0, \ldots , Z_g]$  the homogeneous coordinate ring of $\PP^g$, and consider a finitely generated graded $S$-module $M$. As in \eqref{Resolution}, $M$ has a minimal graded free resolution $E_\bullet$
 \[ \xymatrix{
\ldots \ar[r]  & E_2  \ar[r]  & \ar[r] E_1 \ar[r]  & E_0 \ar[r] &M  \ar[r]&0 ,}
\]
where
$  E_i = E_i(M) = \oplus \, S(-a_{i,j})$.\footnote{We are purposely introducing a  shift in  indexing, so that here our resolutions start in homological degree zero rather than one. The reason for this is that we henceforth wish to view the  resolution \eqref{Resolution} of an ideal $I$ as coming from one of $S/I$ with $E_0 = S$. }  
Write
\[   K_{i,1}(M) \ = \ \big\{
 \text{minimal generators of $E_i(M)$ of degree $i + 1$}
\big  \}.  
\]
This is a finite dimensional vector space whose elements we call $i^{\text{th}}$ syzygies of weight $1$. (The space $K_{i,q}$ of syzygies of weight $q$ are defined analogously.) For instance the ideal $I_C$  of the twisted cubic  $C \subseteq \PP^3$ discussed in \S 1 satisfies
\[   \dim K_{2,1}(I_C) \, = \, 2 \ \ , \ \ \dim K_{1,1}(I_C) \, = \, 3. \]
When $M$ is the $S$-module associated to a coherent sheaf $\FF$ on $\PP^g$, we write simply $K_{i,1}(\FF)$. In particular, the weight one syzygies of of the tangent developable $T$ -- which, as it turns out, govern Conjecture \ref{Folk.Conj} -- are given by $K_{i,1}(\OO_T)$.

 Proposition \ref{OT.as.Kernel} then yields
\begin{corollary} For every $i \ge 1$ one has an exact sequence
\begin{equation} \label{Exact.Seq.K1}
0 \lra K_{i,1}(\OO_T) \lra K_{i,1}(\OO_{\PP^1 \times \PP^1}) \lra K_{i,1}(\Omega_C^1). 
\end{equation}
\end{corollary}
\noi Happily, it is quite easy to work out the two right-hand terms in the exact sequence \eqref{Exact.Seq.K1}.

\vskip 10pt

Let $U$ denote the two-dimensional complex vector space of linear functions on $\PP^1$, so that $\PP^1 = \PP(U)$ is the projective space of one-dimensional quotients of $V$. The group $\SL_2(\CC)$ acts on everything in sight, and in particular the Koszul groups $K_{i,1}$ will be representations of $\SL_2(\CC)$. After choosing an identification $\Lambda^2 U = \CC$, a  standard calculation shows that there is a canonical $\SL_2$-equivariant isomorphism
\begin{equation} \label{Koszul1}   K_{i,1} (\OO_{\PP^1 \times \PP^1} )  \ = \ \Lambda^{i+1}S^{g-2}U \otimes S^{2i}U, \end{equation}
as well as a natural inclusion
\begin{equation} \label{Koszul2}  K_{i,1}( \Omega_C^1) \ \subseteq \ \Lambda^{i+1}S^{g-1}U \otimes S^{i+1}U.\
 \footnote{In arbitrary characteristic, which is the setting considered in \cite{AFPRW}, the computations are more delicate because one has to distinguish between divided and symmetric powers. Working as we are over $\CC$, we can ignore this.} 
 \end{equation}
  In view of Corollary \ref{Exact.Seq.K1}, one then anticipates a mapping
 \begin{equation}\label{Define.gamma}
 \gamma \, : \, \Lambda^{i+1}S^{g-2}U \otimes S^{2i}U \lra \Lambda^{i+1}S^{g-1}U \otimes S^{i+1}U \end{equation}
 whose kernel is $K_{i,1}(\OO_T)$.  AFPRW in effect devote very substantial effort to elucidating  this map, but the upshot is that it is built from several off-the-shelf pieces. To begin with, there is a natural inclusion
 \begin{equation} \label{coWahl}  S^{2i}U \lra \Lambda^2 S^{i+1}U \end{equation}
 which is dual to the so-called Wahl map $\Lambda^2 S^{i+1}U^* \lra S^{2i}U^*$.\footnote{If $W$ is any two-dimensional $C$-vector space with coodinates $x,y$, the Wahl or Gaussian mapping $\Lambda^2 S^{i+1}W \lra S^{2i}W $ is given (up to scaling) by 
$f \wedge g \mapsto \det \left ( 
 \begin{matrix} 
 f_x & f_y \\ g_x & g_y 
 \end{matrix} \right). $}
 Recalling that $\Lambda^{i+1}\big(A \otimes B \big)$ contains $\Lambda^{i+1}A \otimes S^{i+1}B$ as a summand for any vector spaces $A$ and $B$, $\gamma$ then arises as  the composition
\begin{equation}
\begin{gathered}
\xymatrix{
 & \Lambda^{i+1}S^{g-2}U \otimes S^{2i}U \ar[r] &  \Lambda^{i+1}S^{g-2}U \otimes S^{i+1}U \otimes S^{i+1}U   \ar `r[d]  `[l] `[lld]  `[dl] [dl]\\ 
&  \Lambda^{i+1}\big( S^{g-2}U\otimes U \big) \otimes S^{i+1}U \ar[r] & \Lambda^{i+1}S^{g-1}U \otimes S^{i+1}U.
}
 \end{gathered}
 \end{equation}
 We summarize this discussion as
 \begin{theorem} \label{Syzygy.Tang.Summary}
 With $\gamma$ as just specified, $K_{i,1}(\OO_T)$ sits in the exact sequence
 \[ 0 \lra K_{i,1}(\OO_T) \lra \Lambda^{i+1}S^{g-2}U \otimes S^{2i}U \overset{\gamma} \lra \Lambda^{i+1}S^{g-1}U \otimes S^{i+1}U.
 \]
 \end{theorem}
 
 \subsection*{Hermite reciprocity and Koszul modules}
 
Computations  such as \eqref{Koszul1} and \eqref{Koszul2}  are made by studying the cohomology of certain Koszul-type complexes. These can be difficult to deal with because  they involve high wedge powers of a vector space or vector bundle. One of Voisin's key insights was that upon passing to a Hilbert scheme,  complicated multilinear data are encoded into more geometric questions about line bundles. The next step in the proof of AFPRW is an algebraic analogue of this strategy: one uses a  classical theorem of Hermite to reinterpret Theorem \ref{Syzygy.Tang.Summary} in a more tractable  form involving only symmetric products. (In fact the analogy goes farther: a quick proof of Hermite reciprocity proceeds by interpreting $\Lambda^a S^b U$ as the space of global sections of a line bundle on the   
projective space $\PP^a =\Hilb^a(\PP^1) $.)

As above, let $U$ denote a complex vector space of dimension $2$. The result in question is the following.
\begin{reciprocity} For any $a, b > 0$ there is 
 a canonical $\SL_2(\CC)-linear$ isomorphism 
\begin{equation}
\Lambda^a S^b U \ = \ S^{b+1-a} S^a U.
\end{equation}
\end{reciprocity}
\noi (See for example \cite[Exercise 11.35]{FH}.) In positive characteristics this is no longer true, and one of the contributions of \cite{AFPRW} is to give a characteristic-free variant. 

Plugging this into Theorem \ref{Syzygy.Tang.Summary}, one arrives at:
\begin{corollary} \label{Symm.Tang.Syzygies} The Koszul group $K_{i,1}(\OO_T)$ is the kernel of the map
\begin{equation} \label{gamma.prime}
\gamma^\pr :  S^{2i}U \otimes S^{g-i-2}S^{i+1}U \lra S^{i+1}U \otimes S^{g-i-1}S^{i+1}U  \end{equation}
obtained by pulling back the Koszul differential\,\footnote{For any vector space $V$ and integer $a > 0$, there is a natural map
\[ \Lambda^2 V \otimes S^aV \lra V \otimes S^{a+1} V \]
which fits into the longer Koszul-type complex 
\begin{equation} \label{Koszul.footnote}\Lambda^2 V \otimes S^aV \lra V \otimes S^{a+1} V \lra S^{a+2}V
\lra 0. \end{equation}}
\[   \Lambda^{2} S^{i+1}U \otimes S^{g-i-2}S^{i+1}U \lra S^{i+1}U \otimes S^{g-i-1}S^{i+1}U\]
along the ``co-Wahl" mapping $S^{2i}U \lra \Lambda^2 S^{i+1}U$ appearing in \eqref{coWahl}.
\end{corollary}

We now come to one of the main new ideas of \cite{AFPRW}, namely the introduction of Koszul (or Weyman) modules to study \eqref{gamma.prime}. To understand the motivation,  set $V = S^{i+1}U$,  $A = S^{2i}U$, and put $q = g - i - 2$. On the one hand we have from \eqref{coWahl} an inclusion $A \subseteq \Lambda^2 V$, while for $q \ge 0$ there is a Koszul complex
\[ \Lambda^2 V \otimes S^qV \lra V \otimes S^{q+1}V \lra S^{q+2}V.\]
The construction of $\gamma^\pr$ involved splicing these together, giving a three-term complex
\begin{equation} \label{Weyman1}    A \otimes S^qV \overset {\gamma^\pr} \lra V \otimes S^{q+1}V \lra S^{q+2}V \end{equation}
whose left-hand kernel $K = \ker \gamma^\pr$ we would like to understand. Now suppose we knew that \eqref{Weyman1} is exact. Since in any event the map on the right is surjective, this would yield an exact sequence
\[ 0 \lra K \lra  A \otimes S^qV \overset {\gamma^\pr} \lra V \otimes S^{q+1}V \lra S^{q+2}V \lra 0, \]
and we could immediately compute $\dim K$. The very nice observation of AFPRW is that the exactness of \eqref{Weyman1} is essentially automatic provided only that $q \ge \dim V - 3$.  

Turning to details, let $V$ be any complex vector space of dimension $n$, and suppose given a subspace
$ A  \subseteq \Lambda^2 V$. As above, this determines for $q \ge 1$ a three-term complex
\[    A \otimes S^qV \lra V \otimes S^{q+1}V \lra S^{q+2}V  \]
whose homology $W_q(V,A)$ is called (the degree $q$ piece of) the Koszul module associated to $A$ and $V$. The essential result is:
\begin{theorem} \label{Vanishing.Thm.Weyman.Modules}
Assume that no decomposable forms $2$-forms $\eta \in \Lambda^2V^*$ vanish on $A$. Then 
\begin{equation}
W_q(V, A) \ = \ 0 \ \ \text{ for } \ q \, \ge \, \dim V -  3. \end{equation}
\end{theorem}
\noi This was originally proved in characteristic zero in  \cite{AFPRW2} by a relatively painless application of Bott vanishing. An alternative (but largely equivalent) proof in characteristic zero uses vector bundles on projective space and considerations of Castelnuovo--Mumford regularity. In \cite{AFPRW} the argument is extended to positive characteristics.

\begin{remark} \textbf{(Topological applications of Theorem \ref{Vanishing.Thm.Weyman.Modules})}
Before \cite{AFPRW}, the same authors had used Koszul modules in \cite{AFPRW2} to study some interesting topological questions, involving for example K\"ahler groups. \end{remark} 

\subsection*{Completion of the proof} It is now immediate to complete the proof of Folk Conjecture \ref{Folk.Conj}. To begin with, using the symmetry in the resolution of $T$ mentioned following the statement of Green's Conjecture \ref{Green}, one sees that  \ref{Folk.Conj} is equivalent to the assertion that 
\begin{equation} \label{Required.Vanishing}  K_{\left [ \tfrac{g}{2}\right] ,1}(\OO_T) \ = \ 0. \end{equation}
AFPRW treat separately the case of even and odd genus, so suppose that $g = 2n - 3$ is odd. Put 
\[  i \ = \ \left [ \frac{g}{2} \right] \ = \ n-2, \]
 set $V = S^{i+1}U$ -- so that $\dim V = n$ -- and  let $q= g - i - 2 = n - 3$.  Corollary \ref{Symm.Tang.Syzygies}  shows that  $K_{n-2,1}(\OO_T)$ is governed by the complex
\[  S^{2(n-2)}U \otimes S^{n-3}V \lra V \otimes S^{n-2}V \lra S^{n-1}V \]
computing the Weyman module $W_{n-3}(V, S^{2(n-2)}U)$. The hypotheses of Theorem \ref{Vanishing.Thm.Weyman.Modules} are satisfied, and so $W_{n-3}(V, S^{2(n-2)}U) = 0$. Therefore we get an exact sequence
\[ 0 \lra K_{n-2,1}(\OO_T)  \lra  S^{2(n-2)}U  \otimes S^{n-3}V \overset {\gamma^\pr} \lra V \otimes S^{n-2}V \lra S^{n-1}V \lra 0. \]
A computation of dimensions then shows that $\dim K_{n-2,1}(\OO_T) = 0$, and we are done!

 %
 %
 %
 %

 \end{document}